\documentclass[12pt,a4paper]{article}
\usepackage{amsmath}
\usepackage{amsthm}
\usepackage{amssymb}
\usepackage{amsbsy}
\usepackage{amsfonts}
\usepackage{amstext}
\usepackage{amscd}
\usepackage[dvips]{epsfig}

\numberwithin{equation}{section}
\theoremstyle{plain}
\newtheorem{thm}{Theorem}[section]
\newtheorem{prop}[thm]{Proposition}
\newtheorem{cor}[thm]{Corollary}

\theoremstyle{definition}
\newtheorem{exa}[thm]{Example}

\newtheorem{rem}[thm]{Remark}
\newtheorem{defi}[thm]{Definition}

\newcommand{\real}{\mathbb{R}}
\newcommand{\comp}{\mathbb{C}}
\newcommand{\nat}{\mathbb{N}}
\newcommand{\im}{\text{\normalfont Im}}
\newcommand{\re}{\text{\normalfont Re}}
\newcommand{\trr}{\triangleright}

\begin{document}
\title{Fourier and Cauchy-Stieltjes transforms of power laws including stable distributions}

\author{Takahiro Hasebe \\ Graduate School of Science,  Kyoto University,\\  Kyoto 606-8502, Japan\\ E-mail: hsb@kurims.kyoto-u.ac.jp} 
\date{}

\maketitle

\begin{abstract} We introduce a class of probability measures whose densities near infinity are mixtures of Pareto distributions. This class can be characterized by the Fourier transform which has a power series expansion including real powers, not only integer powers. This class includes stable distributions in probability and also non-commutative probability theories. We also characterize the class in terms of the Cauchy-Stieltjes transform and the Voiculescu transform. If the stability index is greater than one, stable distributions in probability theory do not belong to that class, while they do in non-commutative probability.   
\end{abstract}



\section{Introduction} \label{Cauchy}
Probability measures with power-law behaviour have found their applications in a variety of phenomena such as the energy spectrum of fluid~\cite{UF,K1}, the distribution of dark matters~\cite{NFW} and deformed Gaussian distributions in Tsallis statistics~\cite{Tsa1}. 
Readers interested in further information on power laws can find reviews such as \cite{New}. A basic power law in applications is the Pareto distribution whose density is expressed as $cx^{-\alpha-1}$ ($x \geq R >0$) for $\alpha >0$. In particular, a noise is called a $1/f$ noise, or more precisely $1/f^{\alpha+1}$ noise, if the frequency $f$ follows a Pareto distribution with parameter $\alpha$. 
The case $-1 \leq \alpha \leq 0$ also appears in applications. In this case, the integral $\int_R ^\infty f^{-\alpha-1}df$ diverges. In practice, this does not cause a problem since we usually focus on a finite interval $[R, R']$. However, we only consider $\alpha >0$ in this paper for mathematical simplicity. There are a large number of researches on the $1/f$ noise. For instance, an explanation of the origin of $1/f$ noise can be found in \cite{Bak1}.  

In this paper, we introduce a class of probability measures which are mixtures of Pareto distributions. More precisely, let $S$ be an additive sub-semigroup of $[0,\infty)$ such that: 

(S1) $\nat \subset S$; 

(S2) there is a constant $c>0$ such that $\sharp (S \cap [n, n+1)) \leq c^{n+1}$ for any $n \in \nat$.

Then we introduce the following class of probability measures.      

\begin{defi} A probability measure $\mu$ is said to be in $\mathcal{P}_S$ if the following condition is satisfied: 
there exist $\{ a_{\beta}\}_{\beta \in S, \beta>0} \subset \comp$ and $r, R > 0$ such that 
\begin{gather}
\mu|_{|x|\geq R}(dx) = \sum_{\beta \in S, \beta > 0} \im \left(a_\beta \left(\frac{1}{x}\right)^{\beta+1} \right)dx, \label{eq01} \\
|a_\beta| \leq r^\beta,  \label{eq011} \\
r < c^{-1}R, \label{eq012} 
\end{gather}
where $c$ is the constant in the condition (S2). The function $(\frac{1}{x})^{\beta+1}$ is understood to be $e^{i(\beta+1)\pi}|x|^{-\beta -1}$ for $x < -R$.   Under the conditions (\ref{eq011}) and (\ref{eq012}), the series in (\ref{eq01}) is absolutely convergent (see Proposition \ref{prop33}). 
\end{defi}
Let $\mathcal{F}_\mu(z)=\int_{\real}e^{ixz}\mu(dx)$ and $G_\mu(z)=\int_{\real}\frac{1}{z-x}\mu(dx)$ be the Fourier transform and the Cauchy-Stieltjes transform (or Stieltjes transform for short) of a probability measure $\mu$, respectively. 
The main results of this paper are as follows (see Theorems \ref{thm11} and \ref{thm31}). The following are equivalent: 

(1) $\mu \in \mathcal{P}_S$; 

(2) There exist $(c_\gamma)_{\gamma \in S} \subset \comp$ and $A>0$ such that $|c_\gamma| \leq A^\gamma$ for any $\gamma \in S$,  $c_0=1$ and 
\begin{equation}\label{eq00}
\mathcal{F}_\mu(z) = \sum_{\gamma \in S} \frac{c_\gamma}{\Gamma(\gamma+1)} i^\gamma z^\gamma,~ z > 0;  
\end{equation} 

(3) There exist $(d_\gamma)_{\gamma \in S}\subset \comp$ and $A >0$ such that $|d_\gamma| \leq A^\gamma$ for any $\gamma \in S$, $d_0 =1$ and 
\begin{equation}
G_\mu(z) = \sum_{\gamma \in S} \frac{d_\gamma}{z^{\gamma+1}},~z \in \comp_-,~ |z| > cA.
\end{equation}
The above equivalence generalizes a result of \cite{Has5} which corresponds to the case $S = \nat$. The power $z^\gamma$ is defined as the principal value in $\comp \backslash (-\infty,0]$. 
The coefficients $c_\gamma$ and $d_\gamma$ in fact coincide.  
In Subsection \ref{subsec31}, we moreover characterize $\mathcal{P}_S$ in terms of the reciprocal Cauchy transform and the Voiculescu transform. We consider $G_\mu$ on $\comp_-:=\{z \in \comp: \im\, z <0 \}$ rather than on the upper half-plane. This is because the lower half-plane for $G_\mu$ plays the role that the positive real line $(0,\infty)$ does for $\mathcal{F}_\mu$. This correspondence is clarified more in Section \ref{sec3}. 

If we calculate the Fourier transform of the Pareto distribution restricted to $[0,\infty)$, a logarithmic term $z^\beta\log z$ or other singular terms may appear (see Proposition \ref{prop11}). To cancel such a term, we consider Pareto distributions on the negative real line, not only on the positive real line. However, sometimes we do not have to consider the negative part. In Subsection \ref{subsec3}, we clarify when Pareto distributions on the negative real line can be removed, using a Diophantine approximation.

The idea for the equivalence between (1) and (2) comes from stable distributions~\cite{Zol}. 
The distribution $\mu$ of a random variable $X$ is said to be stable if for any $a,b >0$, there are $c > 0$ and $d \in \real$ such that $aX+bX'$ has the same distribution as $c X +d$, where $X'$ is an independent copy of $X$. If $d$ is $0$ for any $a,b$, the distribution $\mu$ is said to be strictly stable. Except for the Gaussian case, the density of a stable distribution shows a power-law tail $cx^{-\alpha-1}$ as $x \to \infty$, $\alpha \in (0,2)$. The parameter $\alpha$ is called an index of stability. A stable distribution has two important aspects: the distribution of a self-similar L\'evy process~\cite{Sat}; the limit distribution of the sum of i.i.d.\ random variables~\cite{K-G}.

A stable distribution $\mu_\alpha$ with stability index $\alpha$ has a Fourier transform of the form 
\begin{equation}\label{eq000}
\mathcal{F}_{\mu_\alpha}(z) = \int_{\real}e^{ixz}\mu_\alpha(dx) =
\begin{cases}\exp \left(i\gamma z -c(1-i\beta \tan(\frac{\alpha \pi}{2})\cdot\text{sign}(z))|z|^\alpha \right), & \alpha \neq 1, \\
\exp \left(i\gamma z -c(1+i\beta \cdot\frac{2}{\pi}(\log|z|)\text{sign}(z))|z| \right), & \alpha = 1, \\
\end{cases}
\end{equation}
where $c \geq 0, \gamma \in \real$ and $\beta \in [-1,1]$~\cite{Sat}.  
We do not consider the case $\alpha =1$ and $\beta \neq 0$. Then we can use a different parametrization: for $z > 0$, 
\begin{equation}\label{eq001}
\mathcal{F}_{\mu_\alpha}(z) = \exp \left(i\gamma z +i^\alpha bz^\alpha \right), \end{equation}
where $b \in \comp \backslash\{0\}$ satisfies $\arg b \in [(1-\alpha)\pi,\pi]$ if $0 < \alpha \leq1$ and $\arg b \in [0,(2-\alpha)\pi]$ if $1 \leq \alpha \leq 2$. It suffices to consider $z>0$ since information on $z <0$ can be recovered from the complex conjugate. 
The Fourier transform $\mathcal{F}_{\mu_\alpha}(z)$ can be expanded in the form 
\begin{equation}\label{eq04}
\mathcal{F}_{\mu_\alpha}(z) = \sum_{\beta \in S_\alpha}^\infty \frac{i^\beta c_\beta}{\Gamma(\beta+1)}z^{\beta},~~~z >0,  
\end{equation}
which is a special case of (\ref{eq00}) with $S_\alpha:=\{m+n\alpha: m,n \in \nat\}$. However, the coefficients $c_\beta$ are bounded as claimed in (2) if and only if $\alpha \in (0,1]$.  

Since the coefficient $c_\gamma~(=d_\gamma)$ in (\ref{eq00}) has an analogy with the moments of a probability measure, we call $c_\gamma$ a $\gamma$-complex moment. This is why we avoid a logarithmic term in the Fourier transform: the concept of $\gamma$-complex moment can be defined thanks to the absence of a logarithmic term.

The idea for the equivalence between (1) and (3) comes from stable distributions in non-commutative probability, or more specifically, in free, monotone and Boolean probability theories~\cite{V1,Mur2,S-W}. 
The reader is referred to Section \ref{sec3} for details and we just mention here a remarkable difference between probability and non-commutative probability theories: for $\alpha \in (1,2]$, the estimate $|c_\beta| \leq A^\beta$ in (\ref{eq00}) does not hold in probability theory, while it does in free, monotone and Boolean probability theories.

\section{The characterization in terms of the Fourier transform}\label{sec2}
\subsection{Generalized power series expansions} 
In this paper, $\nat$ denotes the set of natural numbers $\{0,1,2,3,\cdots\}$. $z^\beta:= e^{\beta \log z}$ for $\beta \in \real$ denotes an analytic function in $\comp \backslash(-\infty,0]$, where $\log z$ is defined so that $-\pi < \im(\log z) < \pi$.  

\begin{rem}\label{rem1}
(1) The following calculations are correct for $z \in \comp \backslash [0,\infty)$: 

(i) $z^\alpha z^\beta = z^{\alpha +\beta}$ for $\alpha, \beta \in \real$, 

(ii) $z^{n\alpha } = (z^\alpha)^n$ for $\alpha \in \real$, $n \in \mathbb{Z}$. 

\noindent
However, $(wz)^\alpha \neq w^\alpha z^\alpha$ and $z^{\alpha\beta} \neq (z^{\alpha})^\beta$ 
in general. 

(2) It is sometimes useful to denote the expansion (\ref{eq04}) as 
\begin{equation}\label{eq02}
\mathcal{F}_\mu(z) = \sum_{m,n \in \nat} \frac{i^{m+n\alpha} c_{m,n}}{\Gamma(m+n\alpha+1)} z^{m+n\alpha},~~~z >0.   
\end{equation}
However, the coefficients $c_{m,n}$ may not be unique. This is because $m_1 + n_1\alpha = m_2 +n_2\alpha$ can occur for distinct $(m_1,n_1)$ and $(m_2,n_2)$. To avoid this non-uniqueness, we introduce summations over a set such as $\sum_{\beta \in S_{\alpha}}\frac{i^\beta c_\beta}{\Gamma(\beta+1)}z^\beta$. Then coefficients $c_\beta$ are unique. 
\end{rem}

\begin{defi}
A series of the form 
\[
\sum_{\gamma \in S}b_\gamma z^{\gamma} ~~~(b_\gamma \in \comp)
\]
 is called a \textit{generalized power series}.  
\end{defi}

We show some examples of the set $S$ satisfying the properties (S1) and (S2). 
\begin{exa}
(1) Let $\alpha >0$ and $S_\alpha:=\{m+n\alpha: m,n \in \nat\}$. Then $S_\alpha$ satisfies (S1) and (S2). The condition (S2) can be checked as follows. If $k \leq m+n\alpha <k+1$, 
then $k-m \leq n\alpha < k-m+1$. For each $m \in \nat$, the number of possible $n$'s is at most $[\frac{1}{\alpha}]+1$. Therefore, $|S_\alpha \cap [k,k+1)| \leq ([\frac{1}{\alpha}]+1)k$.

(2) More generally, let $0 < \alpha_1 < \cdots < \alpha_p$ and $S_{\alpha_1,\cdots,\alpha_p}:=\{n_0+n_1\alpha_1 + \cdots + n_p\alpha_p: n_0,\cdots,n_p \in \nat\}$. Then $S_{\alpha_1,\cdots,\alpha_p}$ satisfies (S1) and (S2): it can be proved that a constant $c >0$ exists such that $|S_{\alpha_1,\cdots,\alpha_p} \cap [k,k+1)| \leq ck^p$. 
\end{exa}

We prove a basic property on convergence of generalized power series. 
   
\begin{prop}\label{prop33}
Let $(a_\gamma)_{\gamma \in S}$ be a sequence of non-negative real numbers. Then the generalized power series 
$\sum_{\gamma \in S} a_\gamma z^\gamma$ converges in $(0, \varepsilon)$ for an $\varepsilon >0$ if and only if there is an $A > 0$ such that $a_\gamma \leq A^{\gamma+1}$.  
\end{prop}
\begin{proof}
We assume that the generalized power series $\sum_{\gamma \in S} a_\gamma z^\gamma$ converges in $(0, \varepsilon)$ with $\varepsilon \in (0,1)$. We note that the estimate 
\[
\begin{split}
\sum_{\gamma \in S} a_\gamma z^\gamma 
&= \sum_{n =0}^\infty \sum_{\gamma \in S\cap[n,n+1)} a_\gamma z^\gamma  
\geq \sum_{n=0}^\infty \left(\sum_{\gamma \in S\cap[n,n+1)} a_\gamma\right) z^{n+1}  
\end{split}
\]
holds. This implies, from the usual theory of series, that a $C  \geq \varepsilon^{-1}$ exists such that $\sum_{\gamma \in S\cap[n,n+1)} a_\gamma \leq C^{n+1}$ for any $n \in \nat$. In particular, $a_\gamma \leq C^{n+1} \leq C^{\gamma+1} $ for $\gamma \in S \cap [n,n+1)$.

Conversely, let us assume that there is an $A > 0$ such that $a_\gamma \leq A^{\gamma+1}$. We can moreover assume that $A >1$. Then 
$\sum_{\gamma \in S} a_\gamma z^\gamma \leq \sum_{n =0}^\infty \sum_{\gamma \in S\cap[n,n+1)} A^{\gamma+1} z^\gamma \leq cA^2\sum_{n =0}^\infty (cAz)^{n}$ is convergent in $(0, (cA)^{-1})$, where $c >0$ is a constant such that $|S \cap [n,n+1)| \leq c^{n+1}$.   
\end{proof}

\subsection{The main result on the Fourier transform} 
Let $[\beta]$ denote the integer part of $\beta$, that is, $[\beta]$ is the largest integer which is not larger than $\beta$. 
First we calculate the Fourier transform of a Pareto distribution. 
\begin{prop}\label{prop11}
(1)  For $\beta \geq 1$, $z>0$ and $R>0$, the following expansion holds. 
\begin{equation}\label{eq11}
\begin{split}
&R^{\beta} \int_R^\infty e^{ixz}x^{-\beta-1}dx \\
&~~~~~~= \sum_{k=1}^{[\beta]-1}\frac{(iRz)^{k-1}}{\beta\cdots(\beta+1-k)}e^{iRz} + \frac{i^{[\beta]-1}c_1(\beta)}{\beta\cdots(\beta-[\beta]+2)}(Rz)^{\beta}\\
&~~~~~~+ \sum_{k\in \nat\backslash\{1,2\}} i^{k+[\beta]-1}\frac{(Rz)^{\beta}-(Rz)^{k+[\beta]-1}}{k!(k-\beta+[\beta]-1)}\cdot\frac{1}{\beta\cdots(\beta-[\beta]+2)} + f_\beta(z), 
\end{split}
\end{equation} 
where $c_1(\beta):=\int_1^\infty e^{ix}x^{-\beta +[\beta]-2}dx$.  $f_\beta(z)$ is defined as follows: 
\[
f_\beta(z)= \left\{
\begin{array}{ll}\displaystyle \frac{1}{\beta\cdots(\beta-[\beta]+2)}\left(i^{[\beta]}\frac{(Rz)^{[\beta]}-(Rz)^{\beta}}{\beta-[\beta]} +\frac{i^{[\beta]+1}}{2}\frac{(Rz)^{[\beta]+1}-(Rz)^{\beta}}{\beta-[\beta]-1}\right),& ~~\beta \notin \nat, \\
\displaystyle \frac{1}{\beta\cdots(\beta-[\beta]+2)}\left(-i^{\beta}(Rz)^{\beta}\log(Rz)+\frac{i^{[\beta]+1}}{2}\frac{(Rz)^{[\beta]+1}-(Rz)^{\beta}}{\beta-[\beta]-1}\right),&~~\beta \in \nat.  
\end{array}
\right. 
\]
If $[\beta]=1$, the first summation of (\ref{eq11}) is understood to be $0$ and $\beta\cdots(\beta-[\beta]+2)$ is understood to be $1$. If $\beta \notin \nat$, $f_\beta(z)$ can be included in the second summation of (\ref{eq11}); however $f_\beta(z)$ is still exceptional for $\beta \notin \nat$ since it may contain a logarithmic term under the limit $\beta \to \infty$. 

(2) For $0 < \beta < 1$, $z>0$ and $R>0$,  the following expansion holds. 
\begin{equation}\label{eq12}
\begin{split}
&R^{\beta} \int_R^\infty e^{ixz}x^{-\beta-1}dx = c_2(\beta) (Rz)^{\beta} + \sum_{k=0}^\infty i^k \frac{(Rz)^k}{k!(\beta-k)}, 
\end{split}
\end{equation} 
where $c_2(\beta):=\int_1 ^\infty e^{ix}x^{-\beta-1}dx + \sum_{k=0}^\infty \frac{i^k}{k!(k-\beta)}$. 
\end{prop}
\begin{proof}
(1) Let us assume that $\beta \geq 2$. For $z>0$ and $R>0$, we have 
\begin{equation}\label{eq13}
R^{\beta}\int_R^\infty e^{izx}x^{-\beta-1}dx = (Rz)^{\beta}\int_{Rz}^\infty e^{ix}x^{-\beta-1}dx. 
\end{equation}\label{eq14}
With integration by parts,  $\int_{Rz}^\infty e^{ix}x^{-\beta-1}dx$ can be computed more. 
\begin{equation}\label{eq15}
\begin{split}
&\int_{Rz}^\infty e^{ix}x^{-\beta-1}dx \\
&~~~~~~~= \frac{(Rz)^{-\beta}}{\beta}e^{iRz} + \frac{i}{\beta} \int_{Rz}^\infty e^{ix}x^{-\beta}dx \\
&~~~~~~~=  \frac{(Rz)^{-\beta}}{\beta}e^{iRz} + \frac{i(Rz)^{-\beta+1}}{\beta(\beta-1)}e^{iRz}+ \frac{i^2}{\beta(\beta-1)} \int_{Rz}^\infty e^{ix}x^{-\beta+1}dx \\ 
&~~~~~~~=\cdots \\
&~~~~~~~= \sum_{k=1}^{[\beta]-1}\frac{i^{k-1}(Rz)^{k-1-\beta}}{\beta\cdots(\beta-k+1)}e^{iRz} + \frac{i^{[\beta]-1}}{\beta\cdots(\beta-[\beta]+2)}\int_{Rz}^\infty e^{ix}x^{-\beta +[\beta]-2}dx. 
\end{split}
\end{equation}
To expand the last integral in terms of powers of $z$, we decompose it into two parts: 
\begin{equation}\label{eq16}
\begin{split}
\int_{Rz}^\infty e^{ix}x^{-\beta +[\beta]-2}dx 
&= \int_1^\infty e^{ix}x^{-\beta +[\beta]-2}dx + \int_{Rz}^1 e^{ix}x^{-\beta +[\beta]-2}dx \\
&= c_1(\beta) + \sum_{k\geq 0,k\neq1} \frac{i^k}{k!}\int_{Rz}^1 x^{k-\beta +[\beta]-2}dx + i\int_{Rz}^1 x^{-\beta + [\beta]-1}dx \\
&= c_1(\beta)+ \sum_{k\geq 0,k\neq 1} \frac{i^k}{k!} \frac{1-(Rz)^{k-\beta +[\beta]-1}}{k-\beta+[\beta]-1} + g_\beta(z), 
\end{split}
\end{equation}
where $g_\beta(z) = i \frac{(Rz)^{-\beta +[\beta]} -1}{\beta-[\beta]}$ if $\beta \notin \nat$ and $g_\beta(z) = -i \log(Rz)$ if $\beta \in \nat$. (\ref{eq11}) follows from the above relations (\ref{eq13})--(\ref{eq16}). 

Let us consider the case $1 \leq \beta < 2$, or equivalently, $[\beta]=1$. Using the same method, one can check that 
\[
R^{\beta}\int_R^\infty e^{ixz}x^{-\beta-1}dx = (Rz)^{\beta}c_1(\beta) + \sum_{k\geq 0, k\neq 1} \frac{i^k}{k!} \frac{(Rz)^{\beta}- (Rk)^{k}}{k-\beta} + i(Rz)^{\beta}\int_{Rz}^\infty x^{-\beta}dx. 
\]
Therefore, the relation (\ref{eq11}) is still true. 

(2) If $0 < \beta < 1$,  the same method is applicable and then (\ref{eq12}) follows. 
\end{proof}

Therefore, a logarithmic term appears in the Fourier transform of $x^{-\beta-1}, x \geq R$ if $\beta \in \nat$. Even in the case $\beta \notin \nat$, a singularity can appear if we consider asymptotic behavior as $\beta \to \infty$. We remove this singularity by taking an appropriate function supported on the negative real line. Then we can characterize the Fourier transforms which can be expanded by $z^\gamma$ for $\gamma \in S$.  

\begin{thm}\label{thm11}
Let $\mu$ be a probability measure. Then $\mu \in \mathcal{P}_S$ if and only if there exist $c_\gamma \in \comp$ and $A >0$ such that:  

(1) $|c_\gamma| \leq A^\gamma$ for any $\gamma \in S, \gamma >0$ and $c_0=1$;  

(2) $\mathcal{F}_\mu(z) = \sum_{\gamma \in S} \frac{c_\gamma}{\Gamma(\gamma+1)} i^\gamma z^\gamma$ for $z > 0$.  
\end{thm}
\begin{rem}
To generalize the expansion $\mathcal{F}_\mu(z) = \sum_{n=0}^\infty\frac{m_n(\mu)}{n!}(iz)^n$ for a compactly supported $\mu$, we use the factor $\frac{1}{\Gamma(\gamma+1)}$ for $\mu \in \mathcal{P}_S$. In Subsection \ref{subsec31}, this generalization turns out to be relevant.   
\end{rem}
\begin{proof}
 We prepare a notation for global behaviour of functions. For complex-valued functions $f,g$ on $S\times \nat \times [0,\infty)$, $f \prec  g$ means that there is a constant $C >0$, independent of $(\beta,n,z)$, such that $|f(\beta,n,z)| \leq C |g(\beta,n,z)|$ for any $(\beta,n,z)$. 

Let us assume that $\mu \in \mathcal{P}_S$. 
By definition, $\mu|_{|x|\geq R}$ can be written as 
\begin{equation}\label{eq17}
\mu|_{|x|\geq R}(dx) = \sum_{\beta \in S, \beta>0} \im \left(a_\beta \left(\frac{1}{x}\right)^{\beta+1} \right)dx,  
\end{equation}
where $|a_\beta| \leq r^\beta$ for an $r \in (0,c^{-1}R)$. 
Let us decompose $\mu$ into three parts: $\mu = \mu_- + \mu_0 + \mu_+$, where $\mu_-:= \mu|_{x \leq -R}$, $\mu_0:=\mu|_{|x|< R}$ and $\mu_+:=\mu|_{x \geq R}$. Since $\mu_0$ is compactly supported, the Fourier transform of $\mu_0$ can be expanded in a series $\sum_{n=0}^\infty \frac{m_n(\mu_0)}{n!}i^nz^n$ with the estimate $|m_n(\mu_0)| \leq R^n$. 

(Step 1) We start from the cancellation of singular terms of $\mathcal{F}_{\mu_+}$ and $\mathcal{F_{\mu_-}}$. Let us calculate the contribution of the Fourier transform of $\mu_-$. From Proposition \ref{prop11}, we observe that  
\[
R^{\beta} \int_{-\infty}^{-R}e^{ixz}x^{-\beta-1}dx = (-1)^{\beta+1} (Rz)^{\beta}\int_{Rz}^\infty e^{-ix} x^{-\beta-1}dx  
\]
with the convention $(-1)^{\beta+1} =e^{i(\beta+1)\pi}$. 
Therefore, $R^{-\beta}\im((-1)^{\beta+1}a_\beta) \overline{f_\beta(z)}$ is the singular term contributed by $\mu_-$. Altogether, the singular part of $\mathcal{F}_\mu$ is $R^{-\beta}\left((\im\, a_\beta) f_\beta(z) + \im((-1)^{\beta+1}a_\beta) \overline{f_\beta(z)}\right)$. 
For $\beta \notin \nat$, we have 
\[
\begin{split}
&(\im\, a_\beta) f_\beta(z) + \im((-1)^{\beta+1}a_\beta) \overline{f_\beta(z)} \\
&~~~~~~= 
\frac{(-i)^{[\beta]}}{\beta\cdots(\beta-[\beta]+2)}  \im\left( a_\beta \cdot\frac{(-1)^{[\beta]} - (-1)^{\beta}}{\beta - [\beta]}\right)\left((Rz)^{[\beta]}-(Rz)^{\beta}\right) \\
&~~~~~~~+ \frac{i^{[\beta]+1}}{\beta\cdots(\beta-[\beta]+2)}\im\left(a_\beta \cdot\frac{1-(-1)^{\beta-[\beta]-1}}{\beta -[\beta]-1}\right)\frac{(Rz)^{[\beta]+1}-(Rz)^{\beta}}{2}.
\end{split}
\] 
If $\beta \in \nat$, the first term is $0$. Therefore, the singular terms of $\mathcal{F}_\mu(z)$ can be bounded as  $ \prec \frac{r^{[\beta]}}{\Gamma(\beta+1)}(z^{[\beta]}+z^{[\beta]+1}) + \frac{r^{\beta}}{\Gamma(\beta+1)}z^{\beta}$.

(Step 2) The estimation of $\mathcal{F}_{\mu_-}$ follows from the complex conjugate of $\mathcal{F}_{\mu_+}$, so that we only have to estimate the coefficient of each $z^\gamma$ appearing in $\sum_{\beta \in S, \beta > 0} r^\beta \int_R^\infty e^{izx}x^{-\beta-1}dx$, 
excluding the singular terms. For each $\beta \geq 1$, first let us focus on 
\[
I_1(\beta):=  \left(\frac{r}{R}\right)^{\beta}\sum_{k=1}^{[\beta]-1}\frac{(iRz)^{k-1}}{\beta\cdots(\beta-k+1)}e^{iRz} - \left(\frac{r}{R}\right)^{\beta}\sum_{k \in \nat \backslash\{1,2\}} \frac{(iRz)^{k+[\beta]-1}}{k!(k-\beta+[\beta]-1)}\cdot\frac{1}{\beta\cdots(\beta-[\beta]+2)} 
\]
which appears in (\ref{eq11}) with additional factor $(\frac{r}{R})^{\beta}$. The coefficient of $z^n$ in $I_1(\beta)$ can be calculated as 
\[
\begin{split}
&Coef(I_1(\beta),z^n):=(iR)^n\left(\frac{r}{R}\right)^{\beta}\sum_{k=0}^{\min\{n,[\beta]-2\}}\frac{1}{\beta\cdots(\beta-k)\cdot(n-k)!} \\
&~~~~~~~~~~~~~~~~~+ (iR)^n\left(\frac{r}{R}\right)^{\beta}\frac{d(\beta)}{\beta\cdots (\beta-[\beta]+2)}\cdot \frac{1}{(n-\beta)\cdot(n-[\beta]+1)!}, 
\end{split}
\]
where
\[
d(\beta)=\begin{cases} 
1,& [\beta] = n+1 \text{~or~} [\beta] \leq n-2,\\
0,&\text{otherwise}. 
\end{cases}
\]
We note that $\frac{1}{\beta\cdots(\beta-k)\cdot(n-k)!} \prec \frac{1}{n!}\cdot \frac{\Gamma(\beta - k) \Gamma(k+1)}{\Gamma(\beta+1)}\cdot\frac{n!}{(n-k)!k!} = \frac{1}{n!}\cdot B(\beta+1,k+1) \cdot \frac{n!}{(n-k)!k!}$, where $B(\beta+1,k+1)$ is the beta function $\int_0^1 x^{\beta}(1-x)^k dx$. The last expression can be bounded by just $\frac{1}{n!}\cdot \frac{n!}{k!(n-k)!}$.  
To estimate the coefficient of $z^n$ in the Fourier transform of $\mu_+$, we have to sum $Coef(I_1(\beta),z^n)$ over $\beta \in S$. 
The summation $\sum_{\beta \in S}(iR)^n\left(\frac{r}{R}\right)^{\beta}\sum_{k=0}^{\min\{n,[\beta]-2\}}\frac{1}{\beta\cdots(\beta-k)\cdot(n-k)!}$ can be estimated as  
\[
\prec \sum_{\beta \in S} \left(\frac{r}{R}\right)^{\beta} \frac{R^n}{n!}\sum_{k=0}^{\min\{n,[\beta]-2\}} \frac{n!}{k!(n-k)!} \prec \frac{(2R)^n}{n!}. 
\]
 Since $\frac{1}{\beta \cdots(\beta-[\beta]+2)}\cdot\frac{1}{(n-\beta)(n-[\beta]+1)!} \prec \frac{1}{n!}\cdot\frac{n!}{(n-[\beta])![\beta]!}$, the summation of the last term of $Coef(I_1(\beta),z^n)$ over $\beta \in S$ 
can be estimated as $\prec \frac{r^n}{n!}$.

(Step 3) To estimate the coefficient of $z^\gamma$, we focus on  
\[
I_2(\beta):=\frac{(r/R)^\beta}{\beta\cdots(\beta-[\beta]+2)}
\left(c_1(\beta)i^{[\beta]-1}(Rz)^{\beta} + 
\sum_{k \in \nat\backslash\{1,2\}}i^{k+[\beta]-1}\frac{(Rz)^{\beta}}{k!(k-\beta+[\beta]-1)}\right), 
\]
which comes from (\ref{eq11}). For $\gamma \in S$, the coefficient of $z^\gamma$ in $I_2(\beta)$ is nonzero only if $\beta =\gamma$ and, in that case, it can be estimated as $\prec \frac{r^\gamma }{\Gamma(\gamma+1)}z^{\gamma}$. 

Through the above steps, the desired generalized power series of $\mathcal{F}_\mu$ has been obtained.

(Step 4) Conversely, let us assume the conditions (1) and (2).  Combining these two, we obtain 
\begin{equation}\label{eq18}
|\mathcal{F}_\mu(z)| \leq \sum_{\gamma \in S} \frac{A^\gamma}{\Gamma(\gamma+1)} z^\gamma \leq Ce^{Bz},~~~~z>0, 
\end{equation}
for constants $C,B >0$. The second inequality is proved as follows. For simplicity, we suppose $A >1$. From the condition (S2), $|S \cap [n, n+1)| \leq c^{n+1}$ for a constant $c>0$ for any $n \in \nat$. If $0 <z < 1$, then $\sum_{\gamma \in S} \frac{A^\gamma}{\Gamma(\gamma+1)} z^\gamma = \sum_{n=0}^\infty\sum_{\gamma \in S\cap[n,n+1)} \frac{A^\gamma}{\Gamma(\gamma+1)} z^\gamma \leq \sum_{n=0}^\infty \frac{(cA)^{n+1}}{\Gamma(n+1)} z^n \leq cAe^{cAz}$. If $z\geq 1$, we similarly get $\sum_{\gamma \in S} \frac{A^\gamma}{\Gamma(\gamma+1)} z^\gamma \leq cAze^{cAz}$.

Now we apply L\'evy's inversion formula to calculate $\mu$, following the proof of \cite{Has5}.  We introduce an analytic function 
\[
L_\mu(z):= \sum_{\gamma \in S} \frac{c_\gamma}{\Gamma(\gamma+1)} z^\gamma
\]
 in $\comp \backslash (-\infty, 0]$. 
 We also define $L_\mu(z)$ on $(-\infty,0]$ to be the limit from the upper half-plane. 

For $x > B$ and $N \in \nat$,  let $f_N^{+}(x) := \frac{1}{2\pi i} \int_0 ^N e^{-xz} L_\mu (z)dz$.  
From the estimate (\ref{eq18}), $f_N^{+}(x)$ converges to $f^{+}(x):= \frac{1}{2\pi i}\sum_{\gamma \in S} \frac{c_\gamma}{x^{\gamma+1}}$ as $N \to \infty$ locally uniformly in $(B, \infty)$. 
By the way, changing the path of the contour integral, one gets $f_N^{+}(x) = \frac{1}{2\pi}\int_0 ^N e^{-ix\xi} \mathcal{F}_\mu (z)dz -  \int_{\Gamma_N} e^{-xz} L_\mu (z)dz$, 
where $\Gamma_N = \{Ne^{i\theta}; 0 \leq \theta \leq \frac{\pi}{2}\}$. The second integral converges to $0$ as $N \to \infty$ locally uniformly. Therefore, $\frac{1}{2\pi}\int_0 ^N e^{-ixz} \mathcal{F}_\mu (z)dz \to \frac{1}{2\pi i}\sum_{\gamma \in S} \frac{c_\gamma}{x^{\gamma+1}}$ locally uniformly in $(B, \infty)$. Summed with the complex conjugate, this convergence implies  
\[
\frac{1}{2\pi}\int_{-N} ^N e^{-ixz} \mathcal{F}_\mu (z)dz  \to  \frac{1}{\pi}\sum_{\gamma \in S} \frac{\im~c_\gamma}{x^{\gamma+1}} 
\]
locally uniformly in $(B,\infty)$. Using $f^{-}_N(x):= \int_{-N}^0 e^{-xz} L_\mu(z)dz$, we can similarly prove  
\[
\frac{1}{2\pi}\int_{-N} ^N e^{-ixz} \mathcal{F}_\mu (z)dz  \to \frac{1}{\pi}\sum_{\gamma \in S} \im\left( c_\gamma \left(\frac{1}{x} \right)^{\gamma+1}\right)
\]
locally uniformly in $(-\infty,-B)$. L\'{e}vy's inversion formula implies that $\mu$ has the absolutely continuous density $\frac{1}{\pi}\sum_{\gamma \in S, \gamma>0} \im\left( c_\gamma \left(\frac{1}{x} \right)^{\gamma+1}\right)$ for $|x| > B$. 
 \end{proof}

\begin{exa}\label{exa13} (1) A Pareto distribution with the density $cx^{-\alpha-1}$ on $x \geq R >0$ belongs to $\mathcal{P}_{S_\alpha}$ for $\alpha > 0, ~\alpha \notin \nat$, where $S_\alpha:=\{m+n\alpha: m,n \in \nat \}$. Indeed, we can take $a_\alpha = \frac{c}{\sin(\alpha+1)\pi}e^{i(\alpha+1)\pi}$ and $a_\beta =0$ for $\beta \neq \alpha$ in (\ref{eq01}). If $\alpha \in \nat$, a logarithmic term appears in the Fourier transform and the Pareto distribution does not belong to $\mathcal{P}_{S_\alpha}$ (see Proposition \ref{prop11} and Theorem \ref{thm11}). In addition, a Pareto distribution is infinitely divisible since its density is completely monotone~\cite{Sat}.  

(2) As mentioned in Introduction, any $\alpha$-stable distribution with $0 <\alpha < 1$ belongs to $\mathcal{P}_{S_\alpha}$ and Cauchy distributions also belong to $\mathcal{P}_{\nat}$. Under the notation of (\ref{eq001}), the density can be written as 
\begin{equation}\label{eq0000}
\frac{1}{\pi}\sum_{(m,n)\in\nat^2\backslash\{(0,0)\}}\gamma^m\frac{\Gamma(m+\alpha n+1)}{m!n!}\im\left(b^n \left(\frac{1}{x}\right)^{m+n\alpha+1}\right),~~~x \neq 0.  
\end{equation}
In the case $1 < \alpha \leq 2$, however, an $\alpha$-stable distribution does not belong to $\mathcal{P}_{S_\alpha}$ since the double series (\ref{eq0000}) is not convergent (the series (\ref{eq0000}) is still true as an asymptotic expansion~\cite{Sat,Zol}). By contrast, $\alpha$-stable distributions in non-commutative probability belong to $\mathcal{P}_{S_\alpha}$ for any $\alpha \in (0,2]$ (see Section \ref{sec3}).   

(3) Let $\nu$ be a probability measure on $[0,\infty)$ and $0 < \alpha \leq 1$. Then the function $\int_0^\infty e^{-|z|^\alpha x} \nu(dx)$ is the Fourier transform of an infinitely divisible distribution $\mu$ (see Corollary 10.6 of \cite{St-Ha}). The probability measure $\mu$ is a mixture of symmetric $\alpha$-stable distributions. If, moreover, $\nu$ is compactly supported, $\mu$ belongs to $\mathcal{P}_{S_\alpha}$. Indeed, we have 
\[
\mathcal{F}_\mu(z) = \sum_{n=0}^\infty\frac{(-1)^n m_n(\nu)}{n!}z^{\alpha n}, ~~~z>0,  
\]
and from L\'evy's inversion formula, the density of $\mu$ is 
\[
\frac{1}{\pi}\sum_{n=1}^\infty (-1)^n m_n(\nu)\frac{\Gamma(\alpha n +1)}{n!} \im\left(e^{-\frac{\alpha n \pi}{2}i}\left(\frac{1}{x}\right)^{\alpha n +1}\right),~~~~x \neq 0. 
\]
This implies that $\mu \in \mathcal{P}_{S_\alpha}$.  
\end{exa} 
We do not treat stable distributions with index one which are not Cauchy distributions. Such a probability measure has a logarithm function in its density, so that it is not in the scope of (\ref{eq01}).

\subsection{Probability measures with supports bounded below}\label{subsec3}
In applications, power laws with supports bounded below form an important class. For instance, some stable distributions are supported on the positive real line. In this subsection, we consider when a probability measure $\mu \in \mathcal{P}_{S_\alpha}$ has a support bounded below. 

Following the definition of $\mathcal{P}_S$, let us introduce a class of probability measures $\mathcal{P}^B_S$ whose supports are bounded below.   
\begin{defi} \label{defi111}
A class $\mathcal{P}^B_S$ of probability measures consists of any probability measure $\mu$  satisfying the following properties: 

(1) the support of $\mu$ is bounded below;  

(2) there exist $\{ b_{\beta}\}_{\beta \in S, \beta>0} \subset \real$ and $r, R > 0$ such that 
\begin{gather}
\mu|_{x \geq R}(dx) = \sum_{\beta \in S, \beta > 0} b_\beta \left(\frac{1}{x}\right)^{\beta+1} dx, \label{eq019} \\
|b_\beta| \leq r^\beta,  \label{eq0119} \\
r < c^{-1}R, \label{eq0129} 
\end{gather}
where $c$ is the constant in the condition (S2). 
\end{defi}
A probability distribution of $\mathcal{P}^B_S$ sometimes contains a logarithmic term in its Fourier transform, like the Pareto distribution with parameter $\alpha \in \nat$. However, a stable distribution on the positive real line does not have a logarithmic term in the Fourier transform. To solve this problem, we focus on the case $S = S_\alpha$ ($\alpha >0$). Then we  consider when $\mu \in \mathcal{P}^B_{S_\alpha}$ is contained in $\mathcal{P}_{S_\alpha}$, which means that the Fourier transform of $\mu$ can be written as (\ref{eq00}). A crucial concept is a Diophantine approximation in number theory, whose importance in the context of stable distributions was clarified in \cite{Kuz1,Kuz2,Kuz3}. 

\begin{defi}\label{defi112}
Let us consider the following condition on $\beta \notin \mathbb{Q}$: for any $b > 1$, 
the inequality 
\begin{equation}\label{eq025}
\left| \beta -\frac{p}{q}\right| < b^{-q}
\end{equation}
holds for infinitely many pairs $(p,q)$ of $\mathbb{Z} \times (\nat \backslash\{0\})$. We denote by $\mathcal{D}$ the set of all $\beta \notin \mathbb{Q}$ satisfying the above condition. 
\end{defi}
\begin{rem}
(1) The Lebesgue measure of $\mathcal{D}$ is zero from Khintchine's theorem~\cite{Lang}. 

(2) The set $\mathcal{D}$ is smaller than the set $\mathcal{L}$ used in \cite{Kuz1,Kuz2,Kuz3}. Indeed, $\mathcal{L}$ is defined as in Definition \ref{defi112}, just by replacing ``for any $b>1$'' by ``for some $b >1$''. 
\end{rem}

The following properties are basic. Let us define $\langle \beta \rangle:= \min\{|\beta -n|: n \in \mathbb{Z}\}$ for $\beta \in \real$. 
\begin{prop}\label{prop012} 
(1) Let $\beta \notin \mathbb{Q}$. Then, $\beta \notin \mathcal{D}$ if and only if there is $A >0$ such that $\frac{1}{|\sin(\pi \beta n)|} \leq A^{n}$ for any $n \in \nat, n \geq 1$. 

(2) If $x \in \mathcal{D}$ and $z \in \mathbb{Q}\backslash\{0\}$, then $zx, z+x \in \mathcal{D}$. 

(3) If $x \in \mathcal{D}$, then $x^{-1} \in \mathcal{D}$. 
\end{prop}
\begin{proof}
(2) can be proved from the definition of $\mathcal{D}$. (3) can be proved in a way similar to \cite{Kuz2}. Let us prove (1). If $\beta \notin \mathcal{D}$, then there are $b>1$ and  $q_0 \in \nat$ such that $|\beta q- p| \geq qb^{-q}$ for $q \geq q_0$ and $p \in \mathbb{Z}$. Therefore, we have $\langle \beta q \rangle \geq qb^{-q}$ for $q \geq q_0$. This implies the existence of $A >0$ such that the inequality $\frac{1}{|\sin(\pi \beta q)|} \leq A^{q}$ holds for any $q \in \nat\backslash\{0\}$. The converse is similarly proved. 
\end{proof}
The main theorem is the following. If the coefficient $b_\beta$ in (\ref{eq019}) is not zero for a $\beta \in \nat$, then the Fourier transform of $\mu$ has a logarithmic term as in Proposition \ref{prop11}. Therefore, we assume that $b_\beta = 0$ for any $\beta \in \nat$.  
\begin{thm} \label{thm1212} Let $\mu \in \mathcal{P}^B_{S_\alpha}$ and $\alpha \notin \mathcal{D}$. If the coefficients $b_\beta$ in (\ref{eq019}) are zero for $\beta \in \nat$, then $\mu \in \mathcal{P}_{S_\alpha}$. Conversely, if $\alpha \in \mathcal{D}$, we can find a probability measure $\mu \in \mathcal{P}^B_{S_\alpha} \cap (\mathcal{P}_{S_\alpha})^c$ with $b_\beta = 0$ for $\beta \in \nat$.  
\end{thm}
\begin{proof}
Suppose that a probability measure $\mu$ satisfies the conditions (\ref{eq019})--(\ref{eq0129}). Then the coefficients $\{a_\beta \}_{\beta \in S_\alpha, \beta >0} \subset \comp$ in (\ref{eq01}) should be written as  
\begin{equation}\label{eq0101}
\im\,a_\beta := b_\beta,~~~~\re\, a_\beta := -\cot(\pi \beta) b_\beta,~~~~\beta \in S_\alpha \backslash \nat.   
\end{equation}
The coefficients $a_\beta$ for $\beta \in \nat$ are defined to be $0$. Let us suppose $\alpha \notin \mathcal{D}$. If $\alpha \in \mathbb{Q}$, then $\sup\{|\cot(\pi \beta)|:\beta \notin \nat \}$ is finite, so that $a_\beta$ can be estimated as $|a_\beta| \leq C^\beta$ for a constant $C >0$. 
If $\alpha \notin \mathbb{Q} \cup \mathcal{D}$, there is a constant $C>0$ such that $|a_\beta| \leq C^\beta$ since $\alpha$ satisfies the estimate in Proposition \ref{prop012}(1). Therefore, in both cases $\alpha \in \mathbb{Q}$ and $\alpha \notin \mathbb{Q}$, $\mu$ can be written as 
\[\mu|_{|x|\geq A}(dx) = \sum_{\beta \in S_\alpha, \beta > 0} \im \left(a_\beta \left(\frac{1}{x}\right)^{\beta+1} \right)dx
\]
for an $A >0$.  

Next, let us assume that $\alpha \in \mathcal{D}$. For instance, let $\mu$
 be defined as 
\[
\mu(dx) = \sum_{\beta \in S_\alpha \backslash \nat} b^\beta x^{-\beta-1} dx
\]
supported on $[R, \infty)$ for some $R >0$ and $b >0$. The coefficients $a_\beta$ in (\ref{eq01}) should be written as (\ref{eq0101}), so that 
\[
\im\,a_\beta := b^\beta,~~~~\re\, a_\beta := -\cot(\pi \beta) b^\beta,~~~~\beta \in S_\alpha \backslash \nat.  
\]
From Proposition \ref{prop012}(1), there is no $A >0$ such that $A^{-\beta-1}|\cot(\pi \beta)|$ is bounded for $\beta \in S_\alpha \backslash \nat$. 
This means that $|a_\beta|$ cannot have an estimate of the form $|a_\beta| \leq D^\beta$, so that $\mu \notin \mathcal{P}_{S_\alpha}$.   
\end{proof}

\begin{exa}
(1) The $\alpha$-stable distribution on $[0,\infty)$ belongs to $\mathcal{P}_{S_\alpha}^B \cap \mathcal{P}_{S_\alpha}$ for any $\alpha \in (0,1)$. If parameters in (\ref{eq000}) satisfy $(\gamma,c,\beta)= (0,\cos(\frac{\pi \alpha}{2}), 1)$, the probability density can be calculated as 
\[
\frac{1}{\pi} \sum_{n=1}^\infty \frac{(-1)^{n-1}\sin(\pi\alpha n)\Gamma(n\alpha +1)}{n!} x^{-1 -n\alpha}  
\]
as a convergent series for $x >0$~\cite{Sat,Zol}. 
 
(2) (The supremum of an $\alpha$-stable process~\cite{Kuz1,Kuz2}) Let $X=(X_t)_{t\geq 0}$ be an $\alpha$-stable process with parameters $(\gamma, c,\beta)$ in (\ref{eq000}) satisfying $\gamma =0$ and $c = \sqrt{1+\beta^2 \tan^2(\frac{\pi \alpha}{2})}$. Let $\rho$ be defined by $\rho = P(X_1 >0)$. The process $S_t:=\sup\{X_s: 0 \leq s \leq t\}$ is called the supremum process. From the self-similarity of $X_t$, the distribution of $S_t$ is the same as that of $t^{1/\alpha}S_1$.  
If $\alpha \in (0,1)$ and $\alpha \notin \mathcal{L}\cup \mathbb{Q}$,  
the probability density of $S_1$ is calculated as a convergent double series  
\[
\frac{d}{dx}P(S_1 \leq x) = x^{-1-\alpha}\sum_{m,n \in \nat} b_{m,n+1} x^{-m-n\alpha}
\]
for $x >0$, where $b_{m,n}$ is defined as 
\[
b_{m,n} = \frac{(-1)^{m+n}}{\Gamma(1+\frac{m}{\alpha}+n) \Gamma (-m-\alpha n)} \prod_{j=1}^m \frac{\sin\left(\frac{\pi}{\alpha}(\alpha \rho + j -1)\right)
}{\sin(\frac{\pi j}{\alpha})}\prod_{j=1}^n \frac{\sin(\pi \alpha(\rho + j -1)}{\sin(\pi \alpha j)}. 
\]
Since $\mathcal{D} \subset \mathcal{L}$, we conclude that the distribution of $S_1$ belongs to $\mathcal{P}^B_{S_\alpha}\cap \mathcal{P}_{S_\alpha}$ from Theorem \ref{thm1212}.

(3) (The last passage time~\cite{Kuz3}) Let $Y=(Y_t)_{t \geq 0}$ be a symmetric $\alpha$-stable process in $\real^d$ whose Fourier transform is $E[e^{iv \cdot Y_t}] = e^{-t||v||^\alpha}$ for $v \in \real^d$, where $||\cdot||$ is the Euclidean norm. We define the last passage time $U_r:=\sup\{t >0: ||Y_t|| < r \}$. Then $U_{ar}$ has the same distribution as $a^\alpha U_r$.  If $1 < \alpha < d$, the probability density of $U_2$ can be calculated as 
\[
\frac{d}{dt}P(U_2 \leq t) = \frac{2}{\alpha \Gamma(\frac{d-\alpha}{2})}\sum_{m \geq 0} \frac{(-1)^m\Gamma(\frac{d+2m}{\alpha})}{m!\Gamma(\frac{d-\alpha}{2}+m+1)}t^{-\frac{d+2m}{\alpha}}, 
\]
which is convergent for $t > 0$. From Proposition \ref{prop012}(3) and Theorem \ref{thm1212}, this distribution belongs to $\mathcal{P}^B_{S_{\alpha^{-1}}}\cap \mathcal{P}_{S_{\alpha^{-1}}}$ if $\alpha \notin \mathcal{D}$. 

More examples can be found in \cite{Kuz3}. 
\end{exa}

\subsection{Generalized moments and basic properties}\label{subsec5}

It is important to mention the uniqueness of the coefficient of $z^\gamma$. 
\begin{prop} \label{prop2} Let $\mu \in \mathcal{P}_S$. If 
$\mathcal{F}_\mu(z) = \sum_{\gamma \in S} \frac{c_\gamma}{\Gamma(\gamma+1)} i^\gamma z^\gamma= \sum_{\gamma \in S} \frac{d_\gamma}{\Gamma(\gamma+1)} i^\gamma z^\gamma$ for $z > 0$, then $c_\gamma = d_\gamma$ for any $\gamma \in S$. 
\end{prop} 
\begin{proof}
This can be proved inductively in terms of asymptotic behavior as $z \searrow  0$. 
\end{proof}

In the paper \cite{Has5}, we introduced the concept of complex moment. The class $\mathcal{P}_S$ enables us to generalize this concept.  
\begin{defi}\label{defi030} Let $S$ be a set satisfying (S1) and (S2). 
The complex number $c_\gamma$ in Theorem \ref{thm11} is called a $\gamma$-complex moment of $\mu$. It is denote by $m_\gamma(\mu)$. If $X$ is a random variable with the distribution $\mu$, we also write $m_\gamma(\mu) = m_\gamma(X)$.  
\end{defi}
\begin{rem}
(1) If $\mu$ is compactly supported and $\gamma \in \nat$, the above $\gamma$-complex moment coincides with the usual moment. Moreover, if $\mu \in \mathcal{P}_{\nat}$, $m_\gamma(\mu)$ coincides with the $\gamma$th complex moment~\cite{Has5}. Therefore, no confusion arises if we use the same symbol $m_\gamma(\mu)$ in this paper. 
\end{rem}
From the proof of Theorem \ref{thm11}, we conclude the following. 
\begin{prop}\label{prop27}
In (\ref{eq01}), we can take $a_\gamma$ to be $a_\gamma = \frac{1}{\pi}m_\gamma(\mu)$. We note that if $\gamma \in \nat$, the real part of $a_\gamma$ is arbitrary. If $\gamma \notin \nat$, $a_\gamma$ is uniquely determined by $\mu$, so that $a_\gamma =\frac{1}{\pi}m_\gamma(\mu)$ is the unique choice. 
\end{prop}

The binomial-type expansion is true for $m_\gamma(\mu)$. 
\begin{prop}\label{prop23} (1) $\mathcal{P}_S$ is closed under the convolution $\ast$. 

(2) For $\mu,\nu \in \mathcal{P}_S$, we have 
\[
m_\beta(\mu \ast \nu) = \sum_{\gamma,\delta \in S: \gamma+\delta=\beta} \frac{\Gamma(\beta+1)}{\Gamma(\gamma+1)\Gamma(\delta+1)} m_\gamma(\mu)m_\delta(\nu).   
\] 
\end{prop}
\begin{proof}
For $\mu, \nu \in \mathcal{P}_S$, let us define $c_\beta:=\sum_{\gamma,\delta \in S: \gamma+\delta=\beta} \frac{\Gamma(\beta+1)}{\Gamma(\gamma+1)\Gamma(\delta+1)} m_\gamma(\mu)m_\delta(\nu)$. Then, formally,  
\begin{equation}\label{eq08}
\sum_{\beta \in S} \frac{c_\beta}{\Gamma(\beta+1)}(iz)^\beta = \mathcal{F}_\mu(z)\mathcal{F}_\nu(z),~~~z>0. 
\end{equation}
To prove that the LHS is an absolutely convergent series, we estimate $c_\beta$ as $|c_\beta| \leq D^\beta$. Let $c>0$ be the constant in the condition (S2), i.e., $|S \cap [n,n+1)| \leq c^{n+1}$, and let $A >0$ be a constant such that $|m_\gamma(\mu)|, |m_\gamma(\nu)| \leq A^\gamma$. Then 
\[
\begin{split}
|c_\beta| 
&\leq \sum_{n=0}^{[\beta]}\sum_{\gamma \in S \cap [n,n+1),~\gamma +\delta=\beta}\frac{\Gamma(\beta+1)}{\Gamma(\gamma+1)\Gamma(\delta+1)}A^{\gamma+\delta} \\
& \leq \sum_{n=0}^{[\beta]}\sum_{\gamma \in S \cap [n,n+1)}\frac{\Gamma(\beta+1)}{\Gamma(\gamma+1)\Gamma(\beta-\gamma+1)}A^{\beta} \\
& \leq \sum_{n=0}^{[\beta]}\sum_{\gamma \in S \cap [n,n+1)}\frac{\Gamma([\beta]+2)}{\Gamma(n+1)\Gamma([\beta]-n+1)}A^{\beta} \\
& \leq ([\beta]+2)\sum_{n=0}^{[\beta]}c^{n+1}\frac{\Gamma([\beta]+1)}{\Gamma(n+1)\Gamma([\beta]-n+1)}A^{\beta} \\
&= c([\beta]+2)(c+1)^{[\beta]}A^{\beta}. 
\end{split}
\]
Therefore, (\ref{eq08}) holds as convergent series. This in addition implies that $\mu \ast \nu \in \mathcal{P}_S$ and $m_\beta(\mu\ast \nu) = c_\beta$. 
\end{proof}

\section{Power laws in non-commutative probability} \label{sec3} 
We focus on free, monotone and Boolean independences which are important independences in  non-commutative probability. In this context, the Stieltjes transform is more relevant as a moment-generating function than the Fourier transform. Therefore, we consider generalized power series expansions for Stieltjes transforms and related transforms.

\subsection{Preliminaries}
We summarize preparatory concepts and results. Free, monotone and Boolean independences for random variables were introduced in \cite{V1}, \cite{Mur2} and \cite{S-W}, respectively. For  probability measures $\mu$ and $\nu$ on $\real$, the free convolution $\mu \boxplus \nu$  is defined as the distribution of $X+Y$ where $X$ and $Y$ are self-adjoint and free independent random variables with distributions $\mu$ and $\nu$, respectively. The monotone convolution $\mu \trr \nu$ and the Boolean convolution $\mu \uplus \nu$ are defined similarly, with free independence replaced by monotone and Boolean independences, respectively.

The function $F_\mu(z)=\frac{1}{G_\mu(z)}$ is called the reciprocal Cauchy transform of a probability measure $\mu$. 
The function $F_\mu$ has a right inverse $F^{-1}_\mu$ defined in $\Gamma_{\eta,M}:= \{z \in \comp_-: \im\, z < -M,~ \im\, z < \eta|\re \,z| \}$ for some $\eta, M >0$. 
The Voiculescu transform $\phi_\mu$ is then defined as $ \phi_\mu(z)=F_\mu^{-1}(z)-z$. The Voiculescu transform in free probability theory plays the role of $\log\mathcal{F}_\mu$ in probability theory. The following characterizations (1), (2) and (3) were proved in \cite{Be-Vo}, \cite{Mur3} and \cite{S-W}, respectively.  
\begin{thm}
For probability measures $\mu$ and $\nu$ on $\real$, we have the following. 

(1) $\phi_{\mu \boxplus \nu} = \phi_\mu + \phi_\nu$ in an open set of the form $\Gamma_{\eta,M}$.  

(2) $F_{\mu \trr \nu} = F_\mu \circ F_\nu$ in $\comp_-$. 

(3) $F_{\mu \uplus \nu}(z) = F_\mu(z) +F_\nu(z) -z$ in $\comp_-$. 
\end{thm}
Therefore, the transforms $\phi_\mu(z)$ and $F_\mu(z) -z$ play the roles that the logarithm of the Fourier transform does in probability theory.  

In free probability, stable distributions were introduced in \cite{Be-Vo} as analogues of stable distributions in probability theory. A probability measure $\mu$ is called $\boxplus$-stable if for any $a,b>0$, there are $c > 0$  and $d \in \real$ such that $(D_a\mu) \boxplus (D_b\mu) = (D_c\mu) \boxplus \delta_d$. If $d$ is $0$ for any $a,b$, then the distribution is said to be strictly $\boxplus$-stable. In the Boolean case, $\uplus$-stable distributions and strictly $\uplus$-stable distributions can be defined by replacing $\boxplus$ with $\uplus$ \cite{S-W,A-H2}. In the monotone case, only strictly stable distributions have been defined~\cite{Has1}. A probability measure $\mu$ is said to be strictly $\trr$-stable if it is $\trr$-infinitely divisible and the corresponding $\trr$-convolution semigroup $\{\mu_t \}_{t \geq 0}$ with $\mu_1 = \mu$ and $\mu_0 = \delta_0$ satisfies the self-similarity: for any $t \geq 0$, there is a $\lambda > 0$ such that $\mu_t = D_\lambda \mu$.\footnote{In \cite{Wang} J.-C.\ Wang defined strictly $\trr$-stable distributions in a different way and proved that its definition is equivalent to the definition of \cite{Has1}.}

The above stable distributions are characterized as follows~\cite{A-H2,Be-Vo,Has1,S-W}. 
\begin{thm} 
(1) If a probability measure $\mu$ is $\boxplus$-stable (resp.\ $\uplus$-stable), then $\phi_\mu(z)$ (resp.\ $z-F_\mu(z)$) is one of the following forms: 

(i) $-\gamma +bz^{1-\alpha}$ in $\comp_-$ where $\gamma \in \real$, $b \in \comp \backslash\{0\}, \arg b \in [(1-\alpha)\pi,\pi]$ and $\alpha \in (0,1)$. 

(ii) $c - b\log z$ in $\comp_-$ where $c \in \comp_+ \cup \real$ and $b \geq 0$. 

(iii) $-\gamma +bz^{1-\alpha}$ in $\comp_-$ where $a \in \real$, $b \in \comp \backslash\{0\},~ \arg b \in [0,(2-\alpha)\pi]$ and $\alpha \in (1,2]$.

(2) Only for strictly $\trr$-stable laws, we change the definitions of powers: $z^\beta$ is defined to be $e^{\beta \log z}$ in $\comp \backslash [0,\infty)$ so that $\im(\log z) \in (-2\pi,0)$. If a probability measure $\mu$ is strictly $\trr$-stable, then $F_\mu(z) = (z^\alpha - b)^{1/\alpha}$ in $\comp_-$ where $(\alpha,b)$ satisfies one of the following conditions: $b \in \comp \backslash\{0\},~ \arg b \in [(1-\alpha)\pi,\pi]$ and $\alpha \in (0,1]$; $b \in \comp \backslash\{0\},~ \arg b \in [0,(2-\alpha)\pi]$ and $\alpha \in [1,2]$; $b=0$. 
\end{thm}
The parameter $\alpha$ is also called a stability index. In the case (ii), a stability index is defined to be one. We do not consider the case $b > 0$ in (ii) as in the case of probability theory to avoid a logarithmic term of the Stieltjes transform. A $1$-strictly stable distribution is a delta measure or a Cauchy distribution in any case of probability theory, free, Boolean and monotone probability theories.

\subsection{Characterizations in terms of the Stieltjes transform and related transforms}\label{subsec31}
To apply the class $\mathcal{P}_S$ to non-commutative probability theory,  
we prove an analogue of Theorem \ref{thm11} for the Stieltjes transform and its reciprocal. 
\begin{thm}\label{thm31}
Let $\mu$ be a probability measure and $c$ be the constant in the condition (S2). Then the following are equivalent. 

(1) $\mu \in \mathcal{P}_S$. 

(2) There are $(d_\gamma)_{\gamma \in S}\subset \comp$ with $d_0=1$ and $A >0$ such that $|d_\gamma| \leq A^\gamma$ for any $\gamma \in S$ and $G_\mu(z) = \sum_{\gamma \in S} \frac{d_\gamma}{z^{\gamma+1}}$ in $\{z \in \comp_-: |z| > cA \}$. 

(3) There exist $(b_\gamma)_{\gamma \in S}\subset \comp$ $(b_0=1)$ and $A >0$ such that 
$|b_\gamma| \leq A^\gamma$ and $F_\mu(z) = z\sum_{\gamma \in S} \frac{b_\gamma}{z^{\gamma}}$ in $\{z \in \comp_-: |z| > cA \}$. 
\end{thm}
\begin{proof}
$(1) \Leftrightarrow (2)$: We assume that $\mu \in \mathcal{P}_S$. From Theorem \ref{thm11} and Definition \ref{defi030}, $\mathcal{F}_\mu(z) $ has the generalized power series expansion $\mathcal{F}_\mu(z) = \sum_{\gamma \in S} \frac{m_\gamma(\mu)}{\Gamma(\gamma+1)} i^\gamma z^\gamma$ for $z > 0$, and moreover we have $|m_\gamma(\mu)| \leq A^\gamma$ for some $A >0$. By the way, we can see that 
\[
\int_0^\infty \mathcal{F}_\mu(z)e^{-yz}dz = -i G_\mu(-iy),~~ y>0. 
\]
The estimate (\ref{eq18}) is applicable to the present case, to conclude that 
\[
\int_0^\infty \mathcal{F}_\mu(z)e^{-yz}dz = \sum_{\gamma \in S}\frac{m_\gamma(\mu) i^\gamma}{y^{\gamma+1}},~~~~y>cA. 
\]
Since $\frac{1}{i}G_\mu(-iy)$ and $\sum_{\gamma \in S}\frac{m_\gamma(\mu) i^\gamma}{y^{\gamma+1}}$ are both analytic in $(cA,\infty)$, we have the equality 
$G_\mu(z) = \sum_{\gamma \in S} \frac{m_\gamma(\mu)}{z^{\gamma+1}}$ in $\{z \in \comp_-: |z| > cA \}$. 

Let us prove the converse statement.  
Since the limit $\lim_{y \searrow 0}G_\mu(x-iy)$ is locally uniformly in $\real\backslash[-cA,cA]$, we can use the Stieltjes inversion formula, to conclude that $\mu$ is absolutely continuous in $\real\backslash [-cA,cA]$ and 
$\mu|_{|x|> cA}(dx) = \frac{1}{\pi}\sum_{\gamma \in S} \im \left(m_\gamma(\mu) \left(\frac{1}{x}\right)^{\gamma+1} \right)dx$. 

(2) $\Leftrightarrow $ (3):  taking the reciprocal of $G_\mu(z)$ in Theorem \ref{thm31}, we obtain 
\[
F_\mu(z) =z\left( 1-\sum_{\gamma>0} \frac{d_\gamma}{z^\gamma} +\left(\sum_{\gamma>0} \frac{d_\gamma}{z^\gamma}\right)^2 - \left(\sum_{\gamma>0} \frac{d_\gamma}{z^\gamma}\right)^3 +\cdots \right). 
\]
It is not easy to estimate the coefficient of $z^\gamma$ in $F_\mu$, so that we consider another proof based on Proposition \ref{prop33}.
We can observe that this series expansion is absolutely convergent for large $|z|$. Indeed, 
if $|z|$ is large enough, the series $1+\sum_{\gamma>0} \frac{|d_\gamma|}{|z|^\gamma} +\left(\sum_{\gamma>0} \frac{|d_\gamma|}{|z|^\gamma}\right)^2 +\cdots$ converges to a finite real number. 
 Therefore, the order of the summands can be changed and the reordered series 
$F_\mu(z) = z\sum_{\gamma \in S} \frac{b_\gamma}{z^\gamma}$ is also absolutely convergent for large $|z|$. From Proposition \ref{prop33}, the coefficients $b_\gamma$ are bounded as $|b_\gamma| \leq A^\gamma$ for an $A >0$. The converse statement is similarly proved. 
\end{proof}

The above result supports the relevance of the factor $\Gamma(\gamma+1)$ in Theorem \ref{thm11}. Indeed, the following is immediate, while it is not trivial a priori. 
\begin{thm} Let $\mu \in \mathcal{P}_S$. If we expand $G_\mu(z) = \sum_{\gamma \in S} \frac{d_\gamma}{z^{\gamma+1}}$ for large $|z|$, $\im~z <0$, then $d_\gamma = m_\gamma(\mu)$. 
\end{thm}
\begin{rem}
(1) This coincidence of $m_\gamma(\mu)$ and $d_\gamma$ supports the definition of $\gamma$-complex moments: one can define the same $\gamma$-complex moments both in terms of Fourier and Stieltjes transforms. 

(2) Theorem \ref{thm11} and Theorem \ref{thm31} explain the reason why $G_\mu$ is defined on $\comp_-$ rather than on the upper half-plane $\comp_+$. If we stated Theorem \ref{thm31} for 
$G_\mu$ on $\comp_+$, then the coefficient $d_\gamma$ would not be equal to $m_\gamma(\mu)$. Therefore, $(0,\infty)$ for the Fourier transform corresponds to $\comp_-$ for the Stieltjes transform. A similar observation is in \cite{Has5} in the case $S=\nat$.  
\end{rem}

In free probability, the Voiculescu transform is crucial to analyze the additive free convolution. Therefore, we characterize $\mathcal{P}_S$ in terms of the Voiculescu transform. 
For general $S$, the inverse function is difficult to treat, so that we only focus on $S_{\alpha_1,\cdots, \alpha_p}$. A key to the proof is to look at the powers $\{z^{\alpha_i}\}_{i=1}^p$ as independent variables, while all of them are functions of $z$. 
\begin{thm} \label{thm45}
Let $0 < \alpha_1 < \cdots < \alpha_p$ for $p \in \nat, p \geq 1$, $S_{\alpha_1,\cdots,\alpha_p}:=\{n_0+n_1\alpha_1 + \cdots + n_p\alpha_p: n_0,\cdots,n_p \in \nat\}$ and $c$ be the constant in the condition (S2). Then the following conditions are equivalent. 

(1) $\mu \in \mathcal{P}_{S_{\alpha_1,\cdots,\alpha_p}}$. 

(2) There exist $(e_\gamma)_{\gamma \in S_{\alpha_1,\cdots,\alpha_p}}\subset \comp$ and $A >0$ such that 
$|e_\gamma| \leq A^{\gamma+1}$ and $\phi_\mu(z) = \sum_{\gamma \in S_{\alpha_1,\cdots,\alpha_p}} \frac{e_\gamma}{z^{\gamma}}$ in $\{z \in \comp_-: |z| > cA \}$. 
\end{thm} 
\begin{proof} 
(1) $\Rightarrow $ (2): let us expand $G_\mu$ as in Theorem \ref{thm31}(2). 
This expansion can also be written as 
\[
G_\mu(z) = \frac{1}{z}\sum_{(n_0, \cdots, n_p) \in \nat^{p+1}} \widetilde{d}_{n_0,\cdots,n_p}\left( \frac{1}{z}\right)^{n_0}\cdots \left(\left(\frac{1}{z}\right)^{\alpha_p} \right)^{n_p},  
\]
where $\widetilde{d}_{n_0,\cdots,n_p}$ may not be unique. In the above, formulae  $z^{\gamma +\beta} = z^{\gamma}z^{\beta}$ and $z^{n\beta} = (z^{\beta})^n$ were used; see Remark \ref{rem1}. Let us define $g(z):=G_\mu(\frac{1}{z})$ for $z \in \comp_+$ with small $|z|$. 
If $z(1+f(z))$ is an inverse of $g(z)$, $f$ satisfies 
\[
 0= f(z) + \sum_{(n_0, \cdots, n_p) \in \nat^{p+1} \backslash\{0\}} \widetilde{d}_{n_0,\cdots,n_p}z^{n_0 + n_1 \alpha_1 + \cdots + n_p \alpha_p}(1+f(z))^{n_0 +n_1 \alpha_1+\cdots \cdots +n_p\alpha_p +1}.   
\]
To prove the existence of such an $f(z)$, we define 
\[
\widetilde{g}(z_0, \cdots, z_p, w):= w + \sum_{(n_0, \cdots, n_p) \in \nat^{p+1} \backslash\{0\}} \widetilde{d}_{n_0,\cdots,n_p}z_0^{n_0}z_1^{n_1} \cdots z_p^{n_p}(1+w)^{n_0 +n_1 \alpha_1+\cdots \cdots +n_p\alpha_p +1},  
\]  
where $(1+w)^\beta$ is defined by the series expansion $\sum_{n=0}^\infty \binom{\beta}{n}w^n$. 
$\widetilde{g}$ is analytic around $0 \in \comp^{p+2}$, and moreover, $\widetilde{g}(0,\cdots,0)=0$ and $\frac{\partial \widetilde{g}}{\partial w}(0,\cdots,0) =1$. 
From implicit function theorem, there is an analytic mapping $\widetilde{f}$ around $0 \in \comp^{p+1}$ such that 
\[
\widetilde{g}(z_0,\cdots,z_p, \widetilde{f}(z_0,\cdots z_p)) = 0.
\]  
$\widetilde{f}$ has an expansion of the form $\widetilde{f}(z_0,\cdots,z_p)= \sum_{(n_0, \cdots, n_p) \in \nat^{p+1} \backslash\{0\}} \widetilde{f}_{n_0,\cdots,n_p}z_0^{n_0} \cdots z_p^{n_p}$ which is convergent around $0$. Then we can define $f$ to be $f(z):=\widetilde{f}(z,z^{\alpha_1},\cdots,z^{\alpha_p})$. 

Thus, the (right) inverse function of $F_\mu$ exists as $\frac{z}{1+f(\frac{1}{z})}$.  The Voiculescu transform $\phi_\mu(z)$ is then equal to $\frac{-zf(\frac{1}{z})}{1+f(\frac{1}{z})}$, which has a generalized power series with the desired form. We note that the expansion of $\frac{1}{1+f(\frac{1}{z})}$ is convergent for large $|z|$ with $z \in \comp_-$ as discussed in the proof of $(2) \Leftrightarrow (3)$ of Theorem \ref{thm31}. 

The converse implication is similarly proved. 
\end{proof}

Now we have basic properties of $\mathcal{P}_S$ or $\mathcal{P}_{S_{\alpha_1,\cdots,\alpha_p}}$ with respect to free, Boolean and monotone convolutions, in addition to Proposition \ref{prop23}.  
\begin{cor}
(1) For any set $S$ with conditions (S1) and (S2), $\mathcal{P}_S$ is closed under $\uplus$ and $\trr$. 

(2) $\mathcal{P}_{S_{\alpha_1,\cdots,\alpha_p}}$ is closed under $\boxplus$.  
\end{cor}
\begin{proof}
The claims for $\uplus$ and $\boxplus$ are immediate from Theorem \ref{thm31}(3) and Theorem \ref{thm45}, respectively. 
The claim for $\trr$ can be proved with an argument similar to that of Theorem \ref{thm31}. For $\mu,\nu \in \mathcal{P}_S$, the function $F_\mu(F_\nu(z))$ has an absolutely convergent series. Therefore, we can reorder the summands and then $F_\mu(F_\nu(z))$ is of the form $z\sum_{\gamma \in S} \frac{h_\gamma}{z^\gamma}$. Since this is absolutely convergent, Proposition \ref{prop33} implies the existence of an $A >0$ such that $|h_\gamma| \leq A^\gamma$. 
\end{proof}

\begin{exa}
(1) Stable distributions for free, Boolean and monotone independences with stability index $\alpha \in (0,1) \cup(1,2]$ belong to $\mathcal{P}_{S_\alpha}$. Cauchy distributions also belong to $\mathcal{P}_\nat$ which are $1$-stable distributions for any independence. 

(2) Probability distributions $\mu_{b,r}^\alpha$ defined by the Stieltjes transforms 
\[
G_{b,r}^{\alpha }(z)=r^{1/\alpha }\left( \frac{1-(1-b(\frac{1}{z})^{\alpha})^{1/r}}{b}\right) ^{1/\alpha}~~~~(z \in \comp_-)   
\]
were investigated in \cite{A-H}. $\mu_{b,r}^\alpha$ is a probability measure if either of the following conditions is satisfied: 

(i) $1 \leq r < \infty$, $0 < \alpha \leq 1$ and $(1-\alpha)\pi \leq \arg b \leq \pi$; 

(ii) $1 \leq r < \infty$, $1 < \alpha \leq 2$ and $0 \leq \arg b \leq (2-\alpha)\pi$. 

\noindent
$\mu_{b,r}^\alpha \in \mathcal{P}_{S_\alpha}$ since $G_{b,r}^{\alpha }(z)$ has a convergent series of the form $\frac{1}{z}\sum_{n=0}^\infty c_n(\alpha, b, r) \left(\frac{1}{z}\right)^{n\alpha}$ for large $|z|$, $z \in \comp_-$. 
\end{exa}

\section*{Acknowledgements} 
The author thanks Izumi Ojima and Hayato Saigo for discussions on stable distributions and physical applications of power laws. He also thanks Hayato Chiba for a question on a logarithmic singularity in the Fourier transform. Uwe Franz brought the author's attention to the relation between Fourier transforms and Stieltjes transforms, which worked successfully in this paper. Octavio Arizmendi's question on power series expansions was useful to improve this paper. Part of this paper was motivated by talks and discussions in 35th Conference on Stochastic Processes and their Applications. In particular, the author learned much from discussions with Juan Carlos Pardo. This work was supported by Grant-in-Aid for JSPS Research Fellows (21-5106).

\end{document}